\documentclass[12pt]{amsart}
\usepackage{amscd,amssymb}
\usepackage[all]{xy}

\title{On subschemes of formal schemes}
\author{Takehiko Yasuda}
\address{Research Institute for Mathematical Sciences, Kyoto University, Kyoto, 606-8502, Japan}

\email{takehiko@kurims.kyoto-u.ac.jp}

\theoremstyle{plain}
\newtheorem{thm}{Theorem}[section]

\newtheorem{prop}[thm]{Proposition}
\newtheorem{cor}[thm]{Corollary}
\newtheorem{lem}[thm]{Lemma}
\newtheorem{prop-defn}[thm]{Proposition-Definition}

\theoremstyle{definition}
\newtheorem{defn}[thm]{Definition}
\newtheorem{expl}[thm]{Example}

\theoremstyle{remark}
\newtheorem{rem}[thm]{Remark}

\def\AA{\mathbb A}
\newcommand{\CC}{\mathbb C}

\newcommand{\ZZ}{\mathbb Z}

\newcommand{\PP}{\mathbb P}

\newcommand{\NN}{\mathbb N}
\newcommand{\No}{\mathbb{N}_0}

\newcommand{\cF}{\mathcal{F}}

\newcommand{\cI}{\mathcal{I}}
\newcommand{\cJ}{\mathcal{J}}

\newcommand{\cL}{\mathcal{L}}

\newcommand{\cO}{\mathcal{O}}

\newcommand{\cU}{\mathcal{U}}

\newcommand{\cX}{\mathcal{X}}
\newcommand{\cY}{\mathcal{Y}}
\newcommand{\cZ}{\mathcal{Z}}
\newcommand{\cW}{\mathcal{W}}

\newcommand{\fa}{\mathfrak{a}}
\newcommand{\fp}{\mathfrak{p}}

\newcommand{\fm}{\mathfrak{m}}
\newcommand{\fn}{\mathfrak{n}}

\newcommand{\Ker}{\mathrm{Ker}\,}

\newcommand{\Spec}{\mathrm{Spec}\,}
\newcommand{\Spf}{\mathrm{Spf}\,}

\newcommand{\pr}{\mathrm{pr}}

\newcommand{\Supp}{\mathrm{Supp}\,}

\newcommand{\red}{\mathrm{red}}

\newcommand{\adic}{\mathrm{adic}}

\begin{document}

\maketitle

\begin{abstract}
We think about what the subscheme of the formal scheme is.
Differently form the ordinary scheme, the formal scheme has
different notions of ``subscheme''. 
We lay a foundation for these notions and compare them.
We also relate them to singularities of foliations. 
\end{abstract}

\section*{Introduction}

In a foundation of the theory of formal schemes,
it is a problem how to define a subscheme of a formal scheme.
Grothendieck \cite{EGA} defined a closed subscheme  of a locally Noetherian formal scheme.
However it is not definitive, and leads to a pathological phenomenon.
The aim of this paper is to compare different notions of ``subschemes'', 
and to complement the theory of formal schemes.
We also relate the pathological phenomenon with singularities of foliations.

We generalize the notion of the closed subscheme of the locally Noetherian formal scheme
defined in \cite{EGA} to arbitrary ambient formal schemes
and define a subscheme as a closed subscheme of an open subscheme. 
This notion has the advantage 
that a subscheme of a Noetherian  formal scheme is Noetherian,
and the disadvantage that 
a subscheme is not generally an open subscheme of a closed subscheme. 
 B.\ Heinzer  first found an example of the last pathological phenomenon \cite[page 1]{AJL-correction}.
We construct a more explicit one (Theorem \ref{thm-example-no-closure}).

As generalizations of subscheme, we define pre-subschemes and pseudo-subschemes.
In particular, with pseudo-subschemes, the pathological phenomenon mentioned above does not occur: 
A pseudo-subscheme of a formal scheme
with some mild condition is
an open subscheme of a closed pseudo-subscheme.
Unfortunately
 a pseudo-subscheme does not inherit the Noetherianity of the ambient formal scheme,
 but inherit instead a variant of Noetherianity, the ind-Noetherianity.
A closed pseudo-subscheme $\cY \hookrightarrow \cX$
coincides with a  pseudo closed immersion of
 Alonso Tarr\'{\i}o,  Jerem\'{\i}as L\'opez
and P\'erez Rodr\'{\i}guez \cite{AJP}, if $\cX$ and $\cY$ are both locally Noetherian. 
Removing this restriction
is essential and inevitable for our aim.
Subschemes,  pre-subschemes,  pseudo-subschemes of  formal 
schemes are all generalizations of  subschemes of  ordinary schemes.
Showing examples, we will see that subschemes, pre-subschemes and pseudo-subschemes
are mutually different.

Formal schemes naturally appear in studies of algebraic foliations, thanks to Miyaoka's formal 
Frobenius theorem \cite{Miyaoka}. 
Jouanolou \cite{Jouanolou} proved that there exist singular algebraic foliations on $\CC^3$
without any formal separatrix at the origin.
Applying this, we construct closed pseudo-subschemes of $\Spf \CC[w][[x,y,z]]$
that are not closed pre-subschemes nor Noetherian.

McQuillan \cite{McQuillan} modified the definition of formal scheme. 
This modification is the right one, and we follow his definition
with additional minor modifications.

\subsection*{Conventions}

We denote by $\NN$ the set of positive integers, 
and by $\No$ the set of non-negative integers.
A {\em ring} means a commutative ring with unit.
An {\em ordinary scheme} means a scheme, distinguished from a formal scheme
that is not a scheme.

\subsection*{Acknowledgement}
I would like to thank Fumiharu Kato for useful discussions.

\section{Preliminaries}

This section contains generalities on formal schemes.
We adopt a slightly different definition of formal schemes
from those in  \cite{EGA} and \cite{McQuillan}.
However most results in this section are found in \cite{EGA} or \cite{McQuillan},
and similar arguments can apply.

\subsection{Admissible  rings}

For a descending chain of ideals of a ring $A$,
\[
 I_1 \supseteq I_2 \supseteq \cdots,
\]
there exists a unique topology on $A$ which makes $A$ a topological ring and 
for which the collection $\{I_i\}_{i\in \NN}$ of ideals is
a basis of (necessarily open) neighborhoods of $0 \in A$. We call this topology the {\em $\{I_i\}$-topology}.
A {\em linearly topologized ring}  is a topological ring 
with the $\{I_i\}$-topology for some descending chain $\{I_i\}_{i \in \NN}$
of ideals.\footnote{This definition is more restrictive than usual. One usually supposes
only that there exists a (possibly uncountable) collection of ideals forming a basis of neighborhoods of $0$.
In most instances, our condition holds, and it makes arguments simpler.

We will put the corresponding assumption on complete modules and formal schemes. 
In \cite{McQuillan}, this condition is not supposed. However some results in op.\ cit, for instance,
Claim 2.6 and Fact 3.3, seem valid only under this condition.
For  McQuillan uses in the proof the fact that if the projective system of
short exact sequences satisfies the Mittag-Leffler condition, then its limit is also exact.
It is true  only if the system is indexed by $\NN$.} 

An open ideal $I $ of a topological ring  is called an {\em ideal of definition}
if every element $f \in I$ is topologically nilpotent (that is, $f^n \to 0$, as $n \to \infty$).\footnote{This 
definition is due to McQuillan \cite{McQuillan}.
The one in \cite{EGA} is more restrictive:
In \textit{op.\ cit.}, an ideal $I$ is an {\em ideal of definition} if
for every open neighborhood  $V$ of $0$, there exists $n \in \NN$ with $I^n \subseteq V$.}\label{footnote-ideal-of-def} 
A linearly topologized ring is called an {\em admissible ring}
if  it is separated and complete, and admits an ideal of definition. 
Every admissible ring has the largest ideal of definition, the ideal of
the topologically nilpotent elements.\footnote{This fails if we adopt the definition in \cite{EGA}.}

A descending chain of ideals of definition in a topological ring
is called a {\em basis of ideals of definition} if it is a basis 
of neighborhoods of $0$.
If $A$ is an admissible ring with the $\{I_i\}$-topology and 
$J \subseteq A$ is an ideal of definition, then the collection  $\{I_i \cap J \}$
of ideals is a basis of ideals of definition.
Every admissible ring thus admits a basis of ideals of definition.

For a ring $A$ and its ideal $I$,
the {\em $I$-adic topology} on $A$ is by definition the $\{I^n\}_{n \in \NN}$-topology.
An admissible ring $A$ is said to be {\em adic} if the topology on $A$
is identical to the $I$-adic topology for some ideal $I \subseteq A$.

\begin{expl}\label{expl-Noetherian-Nonadic}
Let $k$ be a field and let $A:=k[[x,y]]$  be endowed with the 
$\{(xy^n)\}_{n \in \NN}$-topology. 
Then $A$ is admissible and $\{(xy^n)\}_{n \in \NN} $ is a basis of ideals of definition. 
 However $A$ is not adic.
Indeed, for $ i \ge 2$ and $n \ge 0$, $(x^i y^i)$ does not contain $(xy^n)$,
and hence $(x^i y^i)$ is not open.
\end{expl}

\begin{lem}\label{lem-characterization-admissible}
A topological ring $A$ is admissible if and only if 
$A$ is isomorphic to the limit $\varprojlim A_i$
of some projective system of discrete rings,
\[
 A_1 \leftarrow A_2 \leftarrow \cdots,
\]
such that for every $i$, the map $A_{i+1} \to A_i$
is surjective and every element of its kernel is  nilpotent.

\end{lem}

\begin{proof}
It is essentially the same statement as \cite[0, Lem.\ 7.2.2]{EGA}. 
The proof is parallel.
\end{proof}

\begin{cor}\label{cor-limit-open-ideals-same-radical}
Let $A$ be an admissible ring and
\[
 I_1 \supseteq I_2 \supseteq \cdots
\]
a descending chain of open ideals of $A$ (not necessarily a basis of
neighborhoods of $0$). Suppose that for every $i$, $\sqrt{I_i}= \sqrt{I_1}$.
Then   $B:= \varprojlim A/I_i$ is an admissible ring and the natural map $A \to B$
is a continuous homomorphism.  
\end{cor}

\begin{proof}
The projective system of discrete rings
\[
 A/I_1 \leftarrow A/I_2 \leftarrow \cdots
\]
satisfies the condition in Lemma \ref{lem-characterization-admissible}.
Hence $B$ is admissible. For each $i \in \NN$, put 
\[
B \supseteq \hat I _i := \lim_{\substack{\longleftarrow \\ j \ge i}} I_i/I_j .  
\]
Then the ideals $\hat I_i$ form a basis of neighborhoods of $0 \in B$.
The preimage of $\hat I_i$ in $A$ is $I_i$, in particular, open.
It follows that $A \to B$ is continuous. 
\end{proof}

\begin{defn}
An admissible ring $A$ is said to be {\em pro-Noetherian} if 
one of the following equivalent conditions holds:
\begin{enumerate}
\item For every open ideal $I $, the quotient ring $A/I$ is Noetherian.
\item For some basis $\{I_i\}_{i \in \NN}$ of ideals of definition and for every $i \in \NN$,
the quotient ring $A/I_i$ is Noetherian.
\end{enumerate}
\end{defn}

In particular, every Noetherian admissible ring is pro-Noetherian. 

\begin{lem}\label{lem-pro-Noetherian-admissible}
Let $A$ be a pro-Noetherian admissible ring and $I \subseteq A$ an ideal of definition.
Then for any neighborhood $V$ of $0$, there exists $n \in \NN$ with $I^n \subseteq V$.
(Namely, in the sense of \cite{EGA}, $I$ is an ideal of definition 
and $A$ is admissible. See Footnote \ref{footnote-ideal-of-def})
\end{lem}

\begin{proof}
Since $A$ is linearly topologized, we may suppose that $V$ is an ideal.
Then $ A/V $ is Noetherian. Therefore $I(A/V) $ is finitely generated.
Since every element of $I(A/V) $ is nilpotent, so is $I(A/V) $. This means
that for some $n$, $I^n \subseteq V$.
\end{proof}

\begin{lem}\label{lem-adic+pro-Noetherian-Noetherian}
Every pro-Noetherian adic ring $A$ is Noetherian.
Furthermore for every ideal $I$ of definition in $A$,
the topology on $A$ is identical to the $I$-adic topology.
\end{lem}

\begin{proof}
Let $A$ be a pro-Noetherian adic ring and $I \subseteq A$ an ideal such that  $\{I^n\}_{n \in \NN}$
is a basis of ideals of definition. 
By  definition, $A/I$ and $A/I^2$ are  Noetherian. 
Consequently $I/I^2$ is finitely generated, and from \cite[0, Cor.\ 7.2.6]{EGA},
$A$ is Noetherian. 

Let $J$ be an arbitrary ideal of definition.
Then for some $m \in \NN$, $I ^m \subseteq J$. 
Hence for every $n \in \NN$, $I ^{mn} \subseteq J ^n$, and so $J^n$ is open.
Conversely, since $J$ is finitely generated, for every $n \in \NN$, there exists
$m \in \NN$ with $J^m \subseteq I^n$. This proves the lemma.
\end{proof}

\begin{lem}\label{lem-pro-Noetherian}
Let the notations be as in Corollary \ref{cor-limit-open-ideals-same-radical}.
Suppose that $A$ is pro-Noetherian. Then $B$ is also pro-Noetherian.
\end{lem}

\begin{proof}
The admissible ring $B$ has a basis of ideals of definition $\{ \hat I_i\}_i$,
and the quotient rings $B  / \hat  I_i \cong A/I_i$ are Noetherian.
Hence $B$ is pro-Noetherian.
\end{proof}

The following is required in \S 3.

\begin{prop}\label{prop-surjective-adic}
Let $A$ be a Noetherian adic ring and
\[
 I_1 \supseteq I_2 \supseteq \cdots
\]
a descending chain of ideals of definition in $A$ (not necessarily a basis of
ideals of definition),
and let $B := \varprojlim A/I_i$. Suppose that $B$ is adic.
Then the natural map $A \to B$ is surjective. Moreover the topology on $B$
is the $\{ I_1^n B\}$-topology.
\end{prop}

\begin{proof}
From  Lemmas \ref{lem-adic+pro-Noetherian-Noetherian} and \ref{lem-pro-Noetherian},
 $B$ is also a Noetherian adic ring.
Let 
\[
\hat I _i := \lim _{\substack{ \longleftarrow \\ j \ge i }} I_i /I_j ,
\]
and $J := \hat I_1 $.
Since $J$ is an ideal of definition, from Lemma \ref{lem-adic+pro-Noetherian-Noetherian}, 
the topology on $B$ is identical to the $J$-adic topology. 
Since $\{\hat I_i\}_{i \in \NN}$ is also a basis of ideals of definition of $B$, 
for every  $n$, there exists $i \in \NN$ such that $\hat I_i \subset J^{n} $.
Then we have
\[
 B /J^n =( B / \hat I_i ) /  ( J^{n} /   \hat I_i) = (A /I_i)/(I_1/I_i)^{n}.
\]
For $i \gg m \ge n $, the kernel of $B/J^m \to B/J^n$ is 
\[
 (I_1 /I_i)^{n}/(I_1/I_i)^{m} = I_1^{n} (B/J^m)  .
\]
 Besides $B/J = A /I_1$ is clearly  a finitely generated $A$-module.
As a result, the projective systems  $\{A/I_1^{i}\}$ and $\{B/ J^i\}$ satisfy the conditions of 
\cite[0, Prop.\ 7.2.9]{EGA}, and hence  
\[
J^n = \Ker (B \to B / J^n) = I_1^n B.
\]
This shows the second assertion.

The map $ A/I_1  \to B /I_1 B$ is surjective
and $  B $ is separated for the $ \{ I_1^n B\}$-topology.
From \cite[Th.\ 8.4]{Matusmura}, $A \to B$ is surjective.
\end{proof}

\subsection{Formal schemes}

We associate to each admissible ring $A$ a topologically ringed space $\Spf A$, called 
the {\em formal spectrum} of $A$, as follows:
The underlying topological space of $\Spf A$ is the set of {\em open} prime ideals 
and identified with $\Spec A/I$ for every ideal of definition $I$.
The structure sheaf $\cO_{\Spf A}$ of $\Spf A$ is a sheaf of topological rings.
If $\{I_i\}_{i \in \NN}$ is a basis of ideals of definition of $A$,
then $\cO_{\Spf A}$ is defined to be the limit $\varprojlim \cO_{\Spec A/I_i}$ 
of sheaves $\cO_{\Spec A/I_i}$ of pseudo-discrete rings (for sheaves of pseudo-discrete rings,
see \cite{EGA}).

Let $A$ be an admissible ring with $\{I_i\}$ a basis of ideals of definition
 and $x \in A$.
Then we define the {\em complete localization} of $A$ by $x$, denoted $A_{\{x\}}$,
to be the projective limit $\varprojlim  (A /I_i)_x $ of
the localizations $(A/I_i)_x$ of $A/I_i$.
Then $A_{\{x\}}$ is also admissible.
We call $\Spf A_{\{x\}}$ a {\em distinguished open subscheme} of $\Spf A$. 
The distinguished open subschemes of an affine formal scheme
form a basis of open subsets of the underlying topological space.

For an open prime ideal $\fp \subseteq A$, the stalk $\cO_{\Spf A,\fp}$ 
of the structure sheaf $\cO_{\Spf A}$ at $\fp$ is 
\[
 \lim _{\substack{ \longrightarrow \\ x \notin \fp }} A_{\{x\}} .
\]
We can see that it is a local ring as \cite{EGA}.
 Thus the topologically ringed space $\Spf A $ is a locally topologically ringed space.

After \cite{McQuillan}, we call the completion $\hat \cO_{\Spf A, \fp}$ of $\cO_{\Spf A,\fp}$ 
the {\em fine stalk}, which is isomorphic to
\[
 \varprojlim  \cO_{\Spec A/I_i , \fp /I_i} .  
\]
Here $I_i$, $i \in \NN$, form a basis of ideals of definition.
In the \cite{EGA} notation, if we put $S := A \setminus \fp$, then  
\[
   \cO_{\Spf A, \fp} = A _{\{S\}} \text{ and } \hat \cO_{\Spf A, \fp} = A \{S ^{-1}\}. 
\]
From  \cite[0, Prop.\ 7.6.17]{EGA}, the natural map $\cO_{\Spf A, \fp} \to \hat \cO _{\Spf A ,\fp}$
is a local homomorphism. 

An {\em affine formal scheme} is a topologically ringed space that is isomorphic to 
the formal spectrum of an admissible ring.
A {\em formal scheme} is a topologically ringed space $(X,\cO)$ with $X$ the underlying topological space
and $\cO$ the structure sheaf 
such that there exists an open covering $X = \bigcup _{\lambda \in \Lambda}U_\lambda $, and
for every $\lambda \in \Lambda$, the topologically ringed space $(U_\lambda, \cO|_{U_\lambda})$
is an affine formal scheme.\footnote{This definition of formal scheme is different from
those of \cite{EGA} and \cite{McQuillan}, because of the difference of the definition 
of admissible ring. Moreover in \cite{McQuillan}, it is supposed the existence
of basis of subschemes of definition.}
In particular, every formal scheme is locally topologically ringed space.

A {\em morphism} of formal schemes
is a morphism as locally topologically ringed spaces: The maps of stalks
are necessarily local homomorphisms, equivalently the maps of
fine stalks are local homomorphisms.

A continuous homomorphism $A \to B $ of admissible rings 
induces a morphism of formal schemes, $\Spf B \to \Spf A$.
Then the functor $A \mapsto \Spf A$ from the category of admissible rings to
the category of formal schemes is fully faithful.

\subsection{Several Noetherianities}

We now define some finiteness conditions on formal schemes.

\begin{defn}\label{defn-ind-Noetherian}
Let $\cX$ be a formal scheme.
\begin{enumerate}
\item $\cX$ is {\em adic} if every point of $\cX$ admits an affine open
neighborhood  $\Spf A$ with $A$ adic.
\item $\cX$ is {\em quasi-compact (resp.\ top-Noetherian, locally top-Noetherian}  
if the underlying topological space of $\cX$ 
is quasi-compact (resp.\  Noetherian,  locally Noetherian). 
\item $\cX$ is {\em locally pre-Noetherian (resp.\ locally ind-Noetherian, locally Noetherian)} 
if every point of $\cX$ admits an affine neighborhood $\Spf A$ with $A$ Noetherian
(resp.\ pro-Noetherian, Noetherian and adic).
\item $\cX$ is {\em pre-Noetherian (resp.\ ind-Noetherian, Noetherian)} 
if $\cX$ is locally pre-Noetherian (resp.\ locally ind-Noetherian, locally Noetherian) 
and quasi-compact.
\end{enumerate}
\end{defn}

We have the following implications among properties of a formal scheme:
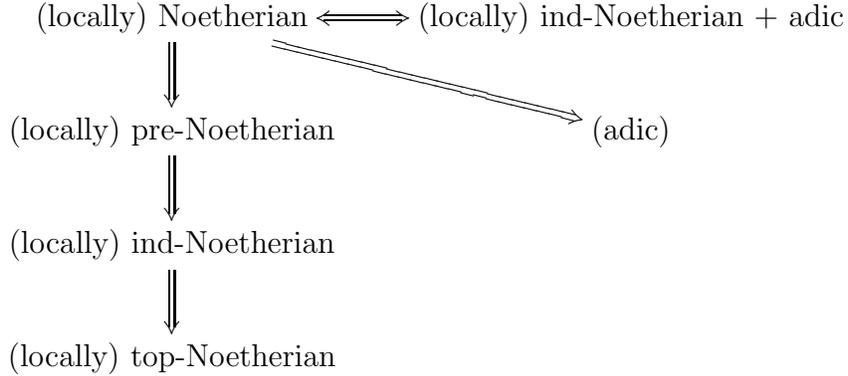
\[\xymatrix{
 \text{(locally) Noetherian} \ar@{=>}[d] \ar@{=>}[dr]  \ar@{<=>}[r] & \text{(locally) ind-Noetherian + adic} \\
 \text{(locally) pre-Noetherian} \ar@{=>}[d]& \text{(adic)} \\
  \text{(locally) ind-Noetherian} \ar@{=>}[d]& \\
   \text{(locally) top-Noetherian}&
}
\]
The top horizontal arrow follows from Lemma \ref{lem-adic+pro-Noetherian-Noetherian}.

\subsection{Quasi-coherent sheaves}

\begin{defn}
Suppose that $A$ is a  linearly topologized ring and
\[
I_1 \supseteq I_2 \supseteq \cdots
\]
a  basis of open ideals.
An $A$-module $M$ endowed with a topology is
said to be {\em complete} if \begin{enumerate}
\item $M$ is a topological group with respect to the given topology and the addition, 
\item there exists a basis of open $A$-submodules of $M$
\[
 M_1 \supseteq M_2 \supseteq \cdots
\]
(that is, the collection $\{M_i\}_{i\in \NN}$ is a basis of open neighborhoods of $0 \in M$)
such that for every $i \in \NN$, $I_i M \subseteq M_i$, and
\item $M$ is separated and complete.
\end{enumerate}
\end{defn}

We note that the second condition is independent of the choice of $\{I_i\}$.

Suppose that $A$ is admissible with a basis $\{I_i\}$ of ideals of definition
and that $M$ is a complete $A$-module
and $\{M_i\}_{i\in \NN}$ is a basis  of open submodules with $I_i M \subseteq M_i$. 
Then $ M /M_i $ is an $A / I_i$-module. 
There exists a corresponding quasi-coherent sheaf $ \widetilde{M /M_i} $ on $\Spec A/I_i$. 
The projective limit 
\[
 M ^\triangle := \varprojlim \widetilde{M /M_i}
\]
of sheaves $ \widetilde{M /M_i} $ of pseudo-discrete groups
is  a  complete $\cO_{\Spf A}$-module, that is, for every open subset $U \subseteq \Spf A$,
$M^\triangle (U)$ is a complete $\cO_{\Spf A}(U)$-module.

\begin{defn}
Let $\cX$ be a formal scheme. 
A  complete $\cO_\cX$-module $\cF$ is said to be {\em quasi-coherent} 
if  every point of $\cX$ has an affine neighborhood $\Spf A$
such that $\cF|_{\Spf A} \cong M^\triangle$ for some complete $A$-module $M$.\footnote{This 
definition is also due to McQuillan (see \cite[\S 5]{McQuillan}), and different
from the \cite{EGA} notion of quasi-coherence.}
\end{defn}

\section{Subschemes of formal schemes and their variants}

In this section, we see various kinds of ``subschemes'' of formal schemes,
compare them, and prove their basic properties.

\subsection{Subschemes}\label{subsec-what-subschemes}

\begin{defn}
An {\em open subscheme} of a formal scheme $\cX$
is an open subset $U$ of $\cX$ along with the restricted structure sheaf $\cO_\cX|_U$.
A morphism $\cY \to \cX$ is said to be an {\em open immersion}
if it is an isomorphism onto an open subscheme of $\cX$.
\end{defn}

\begin{lem}
Let $A$ be an admissible ring with a basis of open ideals $\{J_i\}$ 
and $I \subseteq A$ a closed ideal.
Give the quotient ring $A/I$ the $\{ ( J_i +I)/I \}$-topology.
Then $A/I$ is admissible.
\end{lem}

\begin{proof}
From \cite[the middle of page 56]{Matusmura}, $A/I$ is separated, and
from Theorem 8.1 in op.\ cit, $A/I$ is complete. (Recall that $A$ is supposed to have
a {\em countable} basis of open ideals.)
If $\fa \subseteq A/I$ is the ideal of topologically nilpotent elements,
then its preimage in $A$ contains the largest ideal of definition of $A$,
and open. Therefore $\fa$ is also open (see the top of page 56 in op.\ cit.). 
We conclude that $A/I$ is admissible.
\end{proof}

If $I$ is a closed ideal of an admissible ring $A$, then $I$ is
a complete $A$-module for the induced topology.
Hence we can define an ideal sheaf $I^\triangle \subseteq \cO_{\Spf A}$.
The quotient sheaf $\cO_{\Spf A}/I^\triangle$ is canonically isomorphic to
$\cO_{\Spf A/I}$.

Let $\cX$ be a formal scheme.
An ideal sheaf $\cI \subseteq \cO_\cX$ is said to be {\em closed}
if  every point of $\cX$ has an affine neighborhood
 $\Spf A \subseteq \cX$ such that 
  $\cI|_{\Spf A} = I^\triangle$ for some closed ideal $I \subseteq A$. 
For a closed ideal $\cI \subseteq \cO_\cX$,
the topologically ringed space $ \cY=(\Supp \cO_\cX /\cI , \cO_\cX /\cI )$
is a formal scheme.

\begin{defn}\label{defn-subscheme}
For a formal scheme $\cX$ and a closed ideal sheaf $\cI \subseteq \cO_\cX$,
we call $( \Supp \cO_\cX /\cI , \cO_\cX /\cI  )$ the {\em closed subscheme
defined by $\cI$}. 
A morphism $\cY \to \cX$ of formal schemes is said to be 
a {\em closed immersion}
if it is an isomorphism onto a closed subscheme of $\cX$.
A morphism is said to be an {\em immersion} if it is 
an open immersion followed by a closed immersion.
A {\em subscheme} of a formal scheme $\cX$ is an equivalence class of
immersions $\cY \to \cX$, where $f_i : \cY_i \to \cX$, $i=1,2$, are equivalent
if there exists an isomorphism $g : \cY_1 \to \cY_2$  with $f_1 = f_2 \circ g$.
\end{defn}

Consider the case where $\cX$ is locally Noetherian.
Since every ideal of a Noetherian adic ring is closed 
(see \cite[page 264]{Zariski-Samuel-II}), an ideal sheaf $\cI \subseteq \cO_\cX$
is closed if and only if $\cI$ is coherent in the sense of \cite{EGA}.
Therefore the definition above of the closed subscheme coincides with the one
in \cite{EGA} in this case.

\begin{lem}\label{lem-affine-closed-sub-affine}
Every closed subscheme of an affine formal scheme $\Spf A$ is 
defined by some closed ideal $I \subseteq A$.
\end{lem}

\begin{proof}
Let $\cI \subseteq \cO_{\Spf A}$ be a closed ideal sheaf
and $\{J_i\} _{i \in \NN}$ a basis of open ideals of $A$.
Then for each $i$, there exists an ideal $ I_i \subseteq A / J_i $
such that 
\[
\tilde I_i = \cI + J_i ^\triangle / J_i^\triangle \subseteq \cO_{\Spec A/J_i}.
\]
Moreover the $I_i$ form a projective system and
if we put $I := \varprojlim I_i$, then $I$ is a closed ideal of $A$. 
We easily see that $\cI = I ^\triangle$.
\end{proof}

\begin{prop}
\begin{enumerate}
\item Let $\mathbf{P}$ be any property in Definition \ref{defn-ind-Noetherian} except the quasi-compactness,
and let $\cX$ be a formal scheme satisfying $\mathbf{P}$.
Then every subscheme of $\cX$ satisfies $\mathbf{P}$. 
\item If $g:\cZ \to \cY$ and $f:\cY \to \cX$ are
 (closed) immersions of formal schemes, then $f \circ g$ is also a (closed) immersion.
\item If $\cY \to \cX $ is an immersion of formal schemes 
and $\cZ \to \cX$ is a morphism of formal schemes, then the projection
$\cY \times _\cX \cZ \to \cZ$ is a immersion.
\end{enumerate}
\end{prop}

\begin{proof}
1 and 2 are obvious. 
To prove 3, we may suppose that $\cX = \Spf A$, $\cY = \Spf A/I$ and $\cZ = \Spf B$.
Here $I \subseteq A$ is a closed ideal. Then $\cY \times _\cX \cZ = \Spf ((A /I) \hat \otimes _A B ) $.
Let $\{I_i\}$ and $\{J_i\}$ be bases of ideals of definition in $A$ and $B$ respectively
such that for each $i$, $ I _i$ is contained in the preimage of $J_i$. 
Then we have
\begin{align*}
(A /I) \hat \otimes _A B & \cong \varprojlim  (A/ (I_i +I)  \otimes _{A/I_i} B/J_i) \\
& \cong \varprojlim B/(IB + J_i) \\
& \cong B / \overline {IB} . 
\end{align*}
Here $\overline {IB} $ is the closure of $IB$. Thus $\cY \times _\cX \cZ$ is the closed subscheme
defined by the closed ideal $\overline {IB}$.
\end{proof}

\subsubsection{Pathological examples}

As a consequence of a theorem in \cite{Heinzer-Rotthaus}, 
Bill Heinzer shown the following (see the first page of \cite{AJL-correction}): 

\begin{thm}\label{thm-Heinzer}
 Let $k$ be a field. There exists a nonzero ideal $I \subseteq k[x^\pm,y,z][[t]]$ 
 with $ I \cap  k[x,y,z][[t]]= (0)$.
\end{thm}

The formal scheme $\Spf k[x^\pm ,y,z][[t]]$
is a distinguished open subscheme of $\Spf k[x,y,z][[t]]$
and $\Spf k[x^\pm,y,z][[t]]/I$ is
a subscheme of $\Spf k[x,y,z][[t]]$.
In geometric terms, the theorem  means that the smallest closed subscheme of $\Spf k[x,y,z][[t]]$
containing a subscheme $\Spf  k[x^\pm,y,z][[t]] /I$ of $\Spf k[x,y,z][[t]]$
is $\Spf k[x,y,z][[t]]$ itself. In other words, 
the (scheme-theoretic) closure of  $\Spf  k[x^\pm,y,z][[t]] /I$ in $\Spf k[x,y,z][[t]]$
in a naive sense is $\Spf k[x,y,z][[t]]$.

We find the following simpler and more explicit example:

\begin{thm}\label{thm-example-no-closure}
Consider an element  of  $\CC [x^{\pm },y][[t]]$
\[
 f:= y + a_1 x^{-1}t + a_2 x^{-2}t^2 + a_3 x^{-3}t^3 +\cdots , \ a_i \in \CC \setminus \{0\}.
\]
Suppose that a function $ i\mapsto |a_i| $ is strictly increasing and 
\[
\lim_{i \to \infty} \frac{ |a _{i+1}|}{|a_i|} = \infty.
\] 
Then 
\[
(f) \cap \CC[x,y][[t]] = (0) .
\]
\end{thm}

\begin{proof}
We prove the first assertion by contradiction. So we suppose that
there exists $0 \ne g = \sum _{i \in \No} g_i t^i \in \CC  [x^{\pm },y][[t]]$ with $g_i \in \CC[x^\pm,y]$
such that  $ h := fg \in \CC[x,y][[t]]$. 
If we write $ h = \sum_{i \in \No} h_i t^i $ with $h_i \in \CC[x,y]$,
then for every $i \in \No$, we have
\begin{align*}\label{eq-h-g}
 h_ i  = y g_i + \sum_{j=1}^{i} a_j x ^{- j } g_ {i-j}  = y g_i + h_i' ,\ (h_i'   := \sum_{j=1}^{i} a_j x ^{- j } g_ {i-j}) .
\end{align*}
For each $i \in \No$, write 
\[
 g_i = \sum_{m \in \ZZ, n\in \No} g_{imn} x ^m y^n , \ g_{imn} \in \CC .   
\]
We set 
\begin{align*}
 d_i & := \inf \{ m \in \ZZ | \exists n , g_{imn} \ne 0 \}, \\
 e_i & := \inf \{ n \in \No | g_{id_in} \ne 0  \} , \\
 D_i & := \inf \{ d _ {i-j} -j | 1 \le j \le i   \} \\
 &=\inf \{ d_{j'} -i +j' | 0 \le j' \le i-1 \}, \\
 E_i & := \inf \{ e_{ i-j} | 1 \le j \le i, \ d_{i-j}-j=D_i  \} \\
 &= \inf \{ e_{j'}| 0 \le j' \le i-1, \ 
 d_{j'}-i +j' = D_i \}.
\end{align*}
Here by convention, $\inf \emptyset = + \infty$.
We easily see that for every $i' > i$, 
\[
 D_{i'} < D_i \text{ and } E_{i'} \le E_i.
\]

If for $i_0 \in \NN$, $D_{i_0} < 0$ and 
if  the coefficient of $ x^{D_{i_0}} y^{E_{i_0}} $ in $h'_{i_0}$
is nonzero, then  the coefficient of $ x^{D_{i_0}} y^{E_{i_0}} $ in $yg_{i_0}$
is also nonzero.
Moreover if either ``$m < D_{i_0}$'' or ``$m = D_{i_0}$ and $n < E_{i_0} $'', then
 the coefficient of $ x^m y^n $ in $yg_{i_0}$ vanishes.
It follows that 
\[
  d_{i_0 } = D_{i_0} \text{ and } e_{ i_0} = E_{i_0},
\]
and that
\[
  D_{i_0 +1} = D_{i_0} -1 \text{ and } E_{ i_0+ 1} = E_{i_0}-1,
\]
and that the coefficient of  $ x^{D_{i_0+1}} y^{E_{i_0+1}} $ in $h'_{i_0+1}$
is again nonzero. As a result, $E_{i+1} = E_i -1$ for every $i \ge i_0$.
Since  $E_i \in \No$ for every $i$, it is impossible.

Now it remains to show that for some $i \in \NN$,  $D_i <0$ and 
the coefficient of $x^{D_i}y^{E_i}$ in $h'_i$ is nonzero.
Suppose by contrary that for every $i \in \NN$ with $D_i <0$, 
the coefficient of $x^{D_i}y^{E_i}$ in $h'_i$ is zero. 
Since $i \mapsto D_i$ is strictly decreasing, there exists 
 $ i_1 \in \NN $ such that for every $i \ge i_1$, $D_i < 0$.
Then  for every $i \ge i_1$, the coefficient of $x^{D_i} y ^{E_i}$ in $yg_i $
must be zero.
Therefore we have
\[
 D_ i = D_{i_1} - ( i - i_1 ) \text{ and } E_i = E_{i_1}.
\]
Let  
\[
\Lambda := \{ j| \text{the coefficient of $x^{D_{i_1}} y^{E_{i_1}}$ in $ x^{-j} g_{i_1-j} $ is nonzero} \}
 \subseteq \{ 1,2 ,\dots, i_1 \}
\] 
and let $0 \ne c_j \in \CC $ be the coefficient of $x^{D_{i_1}} y^{E_{i_1}}$ 
in $x^{-j} g_{i_1-j}$, $j \in \Lambda$.
For every $i \ge i_1$, the  coefficient of $x^{D_i}y^{E_i}$ in $h_i'$ is
\[
   \sum _{j \in \Lambda}a_{ j +i-i_1} c_j  = 0.
\]
Let $j_0 \in \Lambda$ be the largest element and $j_1 \in \Lambda$ the second largest one. 
(Note that $\sharp \Lambda \ge 2$).
From the assumption on the  $a_i$, for $i \gg i_1$,
we have
\[
| a_{ j_0 +i-i_1} | - 
( \sharp \Lambda -1)|  a_{ j_1 +i-i_1} | (\max _{j \in \Lambda \setminus \{j_0\} } |c_j /c_{j_0}| )>0.
\]
Therefore, for $i \gg 0$,
\begin{align*}
0&= |\sum _{j \in \Lambda}a_{ j +i-i_1} c_j| \\
& \ge  |c_{j_0}|  \left( | a_{ j_0 +i-i_1} | - \sum_{j \in \Lambda \setminus \{j_0\} }  
|  a_{ j +i-i_1}  c_j  /c_{j_0}  |  \right)\\
& \ge  |c_{j_0}|  \left( | a_{ j_0 +i-i_1} | - 
( \sharp \Lambda -1)|  a_{ j_1 +i-i_1} | (\max _{j \in \Lambda \setminus \{j_0\} } |c_j /c_{j_0}| )\right)\\
& >0 
\end{align*}
This is a contradiction. 
We have proved  the theorem.
\end{proof}

If we remove one more variable, then any ideal as in Theorems \ref{thm-Heinzer} and \ref{thm-example-no-closure} does not exist:

\begin{prop}
Let $k$ be a field.
Then for  any nonzero ideal $I$ of  $ k[x^\pm][[t]]$, 
$ I \cap k[x][[t]] \neq (0)$.
\end{prop}

\begin{proof}
It  suffices to prove the assertion in
 the case where $I$ is principal, say  $ I= (f)$, $f \in k[x^\pm][[t]] $.
Write
\[
f = \sum_{ i \ge n} f_i t^i \in k[x^\pm][[t]], \ f_i \in k[x^\pm], \ f_n \ne 0.
\]
Define  $g_i \in k[x ^\pm, f_n^{-1}]$ inductively as follows;
\[
 g_0 := f_n^{-1} , g _ {i+1} := - ( \sum _{ 0 \le j \le i} g_j f_{n+ i +1 -j} )/ f_n .
\] 
Then 
\begin{align*}
 f   (\sum_{i \ge 0} g_i t^i) &=  \sum _{m \ge n} (( f_n g_{m-n} + \sum _{\substack{i+j =m \\ j < m-n}}f_ig_j ) t^m) \\
 &= t^n + \sum _{m > n} (( - \sum _{j< m-n}g_jf_{m-j} + \sum _{\substack{i+j =m \\ j < m-n}}f_ig_j ) t^m) \\
& = t^n .
\end{align*}
Since $\sum_{i \ge 0} g_i t^i$ is invertible,
ideals $ (f)  $ and $(t^n)$ of $ k[ x^\pm, f_n^{-1}][[t]] $ are identical. 
Glueing $ \Spf k[x^\pm][[t]]/(f) $ and $\Spf k[x, f_n^{-1}][[t]] /(t^n)$, we obtain 
a closed subscheme $\cZ$ of $ \Spf k[x][[t]] $.
Since $ \cZ$ contains $ \Spf k[x^\pm][[t]]/(f)$ as an open subscheme, 
$\cZ$ is not identical to $\Spf k[x][[t]]$.
From Lemma \ref{lem-affine-closed-sub-affine}, it follows that
  $\cZ$ is defined by a nonzero ideal $J \subseteq k[x][[t]]$.
Therefore 
\[
I \cap  k[x][[t]] \supseteq J \neq (0).  
\]
\end{proof}

\subsection{Pre-subschemes}

Recall that a morphism $f : Y \to X$ of ordinary schemes
is a closed immersion if and only if it is a closed embedding
as a map of topological spaces and 
the map $\cO _X \to f_* \cO_Y$ of sheaves
is surjective. We may adopt this as the definition of closed immersions
of ordinary schemes.
However, concerning formal schemes, this condition leads to
a different notion from the closed immersion defined above.

\begin{defn}
A {\em closed pre-immersion} of a formal scheme $\cX$
is a morphism 
$\iota :\cY \to \cX$ such that the map of underlying topological spaces
 is a closed embedding and
 the map $\cO_\cX \to \iota _* \cO_\cY$ is surjective.
An open immersion followed by a closed immersion is 
said to be a {\em pre-immersion}.
 A {\em (closed) pre-subscheme} 
 is an  equivalence class of  (closed) immersions with respect to the equivalence relation
 in Definition  \ref{defn-subscheme}.
 \end{defn}

The following is an example of a closed pre-subscheme that is not a closed subscheme.

\begin{expl}
Let $A := k[[x,y]]$ be endowed with the $\{(xy^i)\}_{i \in \NN}$-topology as in Example \ref{expl-Noetherian-Nonadic} and 
$A^\adic$ the same ring endowed with the $(xy)$-adic topology.
Then $A$ and $A^\adic$ are both admissible rings.
The identity map $A ^\adic \to A$ is a continuous homomorphism.

The formal schemes $ \cX:= \Spf A $ and $\cX ^\adic := \Spf A^\adic $
 have the same underlying topological space,
which consists of three open prime ideals, $(x,y)$,  $ ( x)  $ and $(y)$. 
The stalks of $\cO_\cX$ and $\cO _{\cX^\adic}$ at $ (x,y) $ and $(y)$ are identical as rings,
but not at $ (x) $. 
We have 
\[
\cO_{\cX, (x)} = k((y))[[x]] / (x) =k((y)) \text{ and } \cO_{\cX^\adic, (x)} = k((y))[[x]].
\] 
It follows that the morphism $\cX \to \cX^\adic$
induced by the identity map $A^\adic \to A$ is a closed pre-immersion.
However it is clear that this morphism is not a closed immersion.
\end{expl}

\subsection{Ordinary subschemes}

If $I$ is an open ideal of an admissible ring $A$,
then $I$ is also closed and we can define a formal scheme $\Spf A/I$.
Since $A/I$ is discrete, $\Spf A/I$ is in fact canonically isomorphic to $ \Spec A/I$.
Conversely for a closed ideal $I$ of an admissible ring $A$,
if $\Spf A/I$ is an ordinary scheme, then $I$ is open.

Let $\cX$ be a formal scheme.
An ideal sheaf $\cI \subseteq \cO_\cX$ is said to be {\em open}
if for every point of $\cX$, there exists its affine neighborhood
 $\Spf A \subseteq \cX$ such that $\cI|_{\Spf A} = I^\triangle$
for some open ideal $I \subseteq A$.
Every open ideal sheaf is a closed ideal sheaf.
The closed subscheme of $\cX$ defined by an open ideal is an ordinary scheme. 
 Conversely a closed subscheme that is an ordinary scheme
is defined by an open ideal sheaf.

\begin{defn}
Let $\cX$ be a formal scheme.
A (closed) subscheme $Y \hookrightarrow \cX$ with $Y$  an ordinary scheme
is said to be a {\em (closed) ordinary subscheme}.
\end{defn}

\begin{prop}\label{prop-ordinary-pre-sub}
A pre-subscheme $Y \hookrightarrow \cX$ with $Y$ ordinary subscheme is an ordinary subscheme.
\end{prop}

\begin{proof}
Without loss of generality, we may suppose that $\cX$ is affine, say $\cX = \Spf A$,
and that $Y$ is a closed pre-subscheme. 
Then the underlying topological space of $Y$
is homeomorphic to that of an affine scheme.
Therefore $Y$ is quasi-compact, and covered by  finitely many affine schemes $\Spec B_i$. 
The natural morphism $ \Spec B_i \to \Spf A $ corresponds to a continuous homomorphism
$A \to B_i$. 
Since $B_i$ is discrete, the kernel $ J_i $ of $A \to B_i$ is open. 
Put $ J = \bigcap _i J_i$. Then $J$ is an open ideal and $\iota$ factors as
\[
 Y \xrightarrow{\alpha} \Spec A/J \to \Spf A .
 \]
The $\alpha$ is a closed immersion of ordinary schemes and there exists 
an open ideal $ I \supseteq J $ such that $Y \cong \Spec A/I$. 
Thus $Y$ is a closed ordinary subscheme of $\cX$.
\end{proof}

\subsection{The closure of an ordinary subscheme}\label{subsec-closure-ordinary}

If $Y$ is a subscheme of an ordinary scheme $X$
and if the inclusion map $Y \to X$ is quasi-compact,
then from \cite[Prop.\ 9.5.10]{EGA}, there exists a smallest closed subscheme $\bar Y$ 
 of $X$ that contains $Y$ as an open subscheme. We say that  $\bar Y$
 is the (scheme-theoretic) {\em closure} of $Y$ in $X$.
We can generalize this as follows.

\begin{prop-defn}\label{prop-closure-of-ordinary-subscheme}
Let $\cX$ be a  formal scheme and $Y$ its ordinary subscheme.
Suppose that the inclusion map of the underlying topological spaces
is quasi-compact.
Then there exists a smallest closed ordinary subscheme $\bar Y$ of $\cX$
that contains $Y$ as an open subscheme. 
Moreover, $\bar Y$ is defined by the kernel of $ \cO_\cX \to \iota _* \cO_\cY $,
where $\iota$ is the inclusion. 
We call $ \bar Y$ the {\em closure} of $Y$ in $\cX$.
\end{prop-defn}

\begin{proof}
We first suppose that 
$\cX$ is  affine, say $\cX = \Spf A$.
Then $Y$ is quasi-compact.
Therefore there exists an open subscheme $\cU \subseteq \cX$ covered by
finitely many distinguished open subschemes $\Spf A_{\{f_i\}}$ of $\cX$ such
that $Y$ is a closed ordinary subscheme of $\cU$.
For every $i$, $Y \cap \Spf A_{\{f_i\}}$
is defined by an open ideal $ J_i \subseteq A_{\{f_i\}} $.
If $\{I_ j \}_{j \in \NN}$ is a basis of ideals of definition of $A$,
and if for each $i,j$, 
\[
 (I_j)_{\{f_i\}} := \lim_{\substack{\longleftarrow \\ j' \ge j}} ( I_j /I_{j'})_{f_i}
\]
is the complete localization of $I_j$ by $f_i$,
then $\{ (I_j)_{\{f_i\}} \}_{j \in \NN}$ is a basis of ideals of definition of $A_{\{f_i\}}$.
Hence for $j \gg 0$ and for every $i$, $   (I_j) _{\{f_i\}} \subseteq J_i $.
Then $Y$ is a subscheme of $\Spec A/I_j$.
From \cite[Prop.\ 9.5.10]{EGA}, there exists a
closure $\bar Y$ of $Y$ in $\Spec A/I_j$, which is defined by the
kernel of $\cO_{\Spec A/I_j} \twoheadrightarrow  \iota _* \cO_Y$.
We can view  $\bar Y$  as a closed ordinary subscheme of $\cX $,
which is defined by the kernel of $ \cO_\cX \twoheadrightarrow 
\cO_{\Spec A/I} \twoheadrightarrow  \iota _* \cO_Y $. 
We have proved the assertion in this case.

In the general case, $\cX$ is covered by affine open subschemes $\cX_\lambda$, $\lambda \in \Lambda$. For each $\lambda$, there exists the closure $\bar Y_\lambda$ of
$Y _\lambda := Y \cap \cX_\lambda$ in $\cX_\lambda$. Gluing $\bar Y_\lambda$,
we obtain the closure $\bar Y$ of $Y$ in $\cX$.  
\end{proof}

\subsection{Subschemes of definition}

\begin{defn}
Let $\cX$ be a formal scheme and $\cI \subseteq \cO_\cX$ an open ideal.
We say that $\cI$ is an {\em ideal of definition}
if  every point of $\cX$ has an affine neighborhood $ \Spf A  $ such that
$\cI|_{\Spf A} = I ^\triangle$ for an ideal of definition $I \subseteq A$. 
The ordinary subscheme defined by an ideal of definition
is called a {\em subscheme of definition}.
\end{defn}

\begin{prop}
An ordinary subscheme $Y $ of a formal scheme $\cX$ 
is a subscheme of definition if and only if the underlying topological space
of $Y$ is identical to that of $\cX$.
\end{prop}

\begin{proof}
The ``only if'' direction is trivial.
Suppose that  the underlying topological space
of $Y$ is identical to that of $\cX$.
Without loss of generality, we may suppose, in addition, that 
$\cX$ is affine, say $\cX = \Spf A$.
Let $I \subseteq A$ be the open ideal defining $Y$.
Then  $\Spec A/\sqrt{I}$ is a unique reduced subscheme of $\Spec A$
whose underlying topological  space is identical to that of $\Spf A$.
This shows that $\sqrt{I}$ must be the largest ideal of definition.
Therefore $I$ consists of topologically nilpotent elements, and is an ideal of definition.
\end{proof}

Since every admissible ring $A$ admits a largest ideal of definition,
every affine formal scheme admits  a smallest subscheme of definition, which
is the reduced subscheme of definition.
Glueing the smallest subschemes of definition of affine open subschemes,
we obtain a smallest subscheme of definition of an arbitrary formal scheme.
In particular, every formal scheme has at least one subscheme of definition.

\begin{defn}
Let $\cX$ be a formal scheme and 
\[
 \cI_1 \supseteq \cI_2 \supseteq \cdots
\]
a descending chain of ideals of definition in $\cO_\cX$.
We say that $\{\cI_i\}_{i \in \NN}$ is a {\em basis of ideals of definition}
if every point of $\cX$ has an affine neighborhood $\Spf A $
and there exists a basis $\{I_i\}$ of ideals of definition in $A$ such that for every $i$, $\cI_i|_{\Spf A}= I_i^\triangle$.
The ascending chain of subschemes of definition
corresponding a basis of ideals of definition is called a {\em basis of subschemes of definition}.
\end{defn}

\begin{prop}\label{prop-existence-subscheme-of-definition}
Every top-Noetherian formal scheme has  a basis of subschemes of definition.
\end{prop}

\begin{proof}
Let $\cX$ be a top-Noetherian formal scheme and $\cX = \bigcup _{i=1}^n \cU_i$
its finite affine covering. 
For each  $i$, there exists a basis  $\{\cI_{ij}\}_{j \in \NN}$
 of ideals of definition  on $\cU_i$ and the corresponding basis $\{Y_{ij}\}_{j \in \NN}$
 of subschemes of definition.
For each $i,j$, we denote by $\bar Y_{ij}$ to be the closure of $Y_{ij}$ in $\cX$
and by $\bar \cI_{ij} \subseteq \cO_\cX$ the corresponding ideal sheaf.
For each $j \in \NN$, set $ \cJ_j := \bigcap _{i=1}^n \bar \cI_{ij}$.
The $\cJ_j$'s are  open ideals.
For each $i$,  $\cJ_j|_{\cU_i}$ is contained in $\cI_{ij}$. It follows that 
for each $i$,
$\{\cJ_j|_{\cU_i}\}_{j\in \NN}$ 
is a basis of  ideals of definition and 
so is $\{\cJ_i\}_{i \in \NN}$.
\end{proof}

\begin{prop}
Every locally Noetherian formal scheme has a basis of subschemes of definition.
\end{prop}

\begin{proof}
Let $\cX$ be a locally Noetherian formal scheme and $\cI \subseteq \cO_\cX$ 
the largest ideal of definition.
Then $\{\cI^ n \}_{n \in \NN}$ is a basis of ideals of definition.
\end{proof}

\subsection{Pseudo-subschemes}

\begin{defn}
A {\em pseudo-immersion} of a formal scheme $\cX$
is a morphism $\iota : \cY \to \cX$ of formal schemes 
such that for every immersion $Y \hookrightarrow \cY$ with $Y$ ordinary scheme, 
the composition $Y \hookrightarrow \cY \xrightarrow{\iota} \cX$ is an immersion.
A {\em pseudo-subscheme} is an equivalence class of pseudo-immersions.
\end{defn}

If $\cY$ is a pseudo-subscheme of $\cX$ and if $Y$ is a subscheme of definition,
then $Y$ is by definition an ordinary subscheme of $\cX$.
Therefore the underlying topological space of $\cY$
is a locally closed subset of that of $\cX$.

If $\cY \to \cX$ is a pre-immersion, then
from  Proposition \ref{prop-ordinary-pre-sub},
for every immersion  $Z \hookrightarrow \cY$ with $Z$ ordinary scheme,
the composition $Z \to \cY \to \cX$ is  an immersion.
Hence $\cY \to \cX$ is a pseudo-immersion.

When $\cY$ admits a basis $\{Y_i\}_{i \in \NN}$
of subschemes of definition, then $\cY \to \cX$ is a pseudo-immersion of $\cX$
if and only if for every $i$, $Y_i \hookrightarrow \cY \to \cX$ is an immersion.

\begin{defn}
A pseudo-immersion or a pseudo-subscheme is said to be {\em closed}
if the map of underlying topological spaces is a homeomorphism
onto a closed subset.
\end{defn}

\begin{expl}
Let $ X $ be an ordinary
scheme and $Y$ its closed subscheme. Then the completion $ X_{/Y}$ of $X$ along $Y$ 
is a closed pseudo-subscheme of $X$. 
\end{expl}

Let 
\[
 Z_1 \subseteq Z_2 \subseteq \cdots 
\]
be ordinary subschemes
of a  formal scheme $\cX$, all of which have the same underlying topological space.
Then, from Corollary \ref{cor-limit-open-ideals-same-radical}, the inductive limit
\[
\cZ:= \varinjlim Z_i
\]
is a formal scheme and a pseudo-subscheme of $\cX$.

\begin{expl}\label{expl-infinite-embedded-points}
Suppose that $k$ is an algebraically closed field and 
that a ring $k[x][[t]]$ is endowed with the $(t)$-adic topology.
Let $\cX := \Spf k[x][[t]]$. 
The underlying topological space of $\cX$ is identified with that of $\AA^1_k = \Spec k[x]$.
For each $a \in k$,
we define a subscheme of definition of $\cX$,
\[
Y_a := \Spec k[x][[t]]/  (t^2,(x-a)t).   
\]
It has an embedded point at a rational point $a \in \AA_k^1$. 
For a finite subset $\{a_1,\dots,a_n\} $ of $k$, 
we define  $Y_{a_1,\dots,a_n}$
to be the subscheme of definition of $\cX$ that is
isomorphic to $Y_{a_i}$ around $a_i$, $1 \le i \le n$,
and to $\AA_k^1$ around any point other than $a_1,\dots,a_n$.

Let $\{ a_i ;i\in \NN\}$ be a countable subset of $k$.
 Then we have an ascending chain of subschemes of $\cX$,
\[
 Y_{a_1} \subseteq Y_{a_1,a_2} \subseteq Y_{a_1,a_2,a_3} \subseteq \cdots ,
\]
and obtain a closed pseudo-subscheme of $\cX$, 
\[
\cY :=\varinjlim Y_{a_1,a_2 ,\dots,a_n}.
\]
Then 
\begin{align*}
\hat \cO_{\cY,p} \cong 
\begin{cases}
( k[x,y]/(y^2,xy))_{(x,y)} & (p \in  \{a_1,a_2,\dots\}) \\
k[x]_{(x)} & (p \in \AA_k^1(k) \setminus \{a_1,a_2,\dots\}) \\
k(x) & (p \text{ the generic point}).
\end{cases}
\end{align*}
Thus all fine stalks of $\cO_\cY$ are discrete. 
If $\cY$ is Noetherian, then it is impossible that infinitely many fine stalks of $\cO_{\cY}$
have an embedded prime.
Therefore  $\cY$ is not a closed subscheme of $\cX$. 
Moreover for every open subscheme $\cU \subseteq \cX$, 
$\cY \cap \cU$ is not a closed subscheme of $\cU$ either.
\end{expl}

\begin{prop}\label{prop-properties-pseudo-sub}
\begin{enumerate}
\item
Every pseudo-subscheme of a (locally) ind-Noetherian formal scheme
is (locally) ind-Noetherian.
\item If $g: \cZ \to \cY $ and $f:\cY \to \cX$
are pseudo-immersions, then $f \circ g$ is also a pseudo-immersion.
\item Let $\cY \to \cX$ be a pseudo-immersion of formal schemes and
$\cW \to \cX$ a morphism of formal schemes. Then the projection 
$ \cY \times _\cX \cW \to \cW $ is a  pseudo-immersion.
\item Every closed pseudo-subscheme of an affine formal scheme
is an affine formal scheme.
\end{enumerate}
\end{prop}

\begin{proof}
1. The assertion follows from Proposition \ref{prop-surjective-adic}.

2. If $W \hookrightarrow \cZ$ is an immersion with $W$ ordinary scheme,
then the natural morphism $W \to \cX$ is an immersion.
Therefore $\cZ \to \cX$ is a pseudo-immersion.

3. Since the problem is local, we may suppose that 
$\cY$, $\cX$ and $\cW$ are affine. 
Let $\{Y_i\}$, $\{X_i\}$ and $\{W_i\}$
be bases of subschemes of definition of $\cY$, $\cX$ and $\cW$ respectively
such that for every $i$, the natural morphisms $Y_i \to \cX$ and $W_i \to \cX$
factors through $X_i$. 
Then $ \{Y_i \times _{X_i} W_i\} $ is a basis of subschemes of definition of $\cY \times _\cX \cW$.
Since $Y_i$ is a subscheme of $X_i$, the projection 
$Y_i \times _{X_i} W_i \to W_i$ is an immersion \cite[Prop.\ 4.4.1]{EGA}, 
the natural morphism $ Y_i \times _{X_i} W_i \to \cW $ is also an immersion. 
Therefore $\cY \times _\cX \cW \to \cW$ is a pseudo-immersion.

4. Let $\cX = \Spf A$ and $\cY$ its closed pseudo-subscheme.
Then the underlying topological space of $\cY$
is isomorphic to that of an affine scheme $\Spec R$.
There exists an open covering $\cY = \bigcup _{j=1}^n  \cY_j$
such that for each $j$,  $\cY_j $ is identified with
a distinguished open subscheme $\Spec R_f$ as a topological space. 
Since the $\cY_j \hookrightarrow \cY$ are quasi-compact, 
as in the proof of Proposition \ref{prop-existence-subscheme-of-definition},
we can show that $\cY$ has a basis of subschemes of definition, $\{Y_i\}_{i \in \NN}$.
Then $Y_i$ can be viewed as a closed ordinary subscheme of $\Spf A$.
From Lemma \ref{lem-affine-closed-sub-affine}, the $Y_i$ are affine, say $Y_i = \Spec A_i$.
It follows that $\cY = \Spf (\varprojlim A_i)$ is affine.
\end{proof}

\subsection{Chevalley's theorem}

\begin{thm}\label{thm-single-point}
Let $\cX$ be a Noetherian formal scheme. 
Every pseudo-subscheme of $\cX$ is a subscheme of $\cX$
if and only if the underlying topological space of $\cX$ is discrete.
\end{thm} 

\begin{proof}
The ``if'' direction is essentially due to Chevalley \cite{Chevalley}.
To show this, we may suppose that the underlying topological space of $\cX$ consists of a single point.
Then for some Noetherian complete local ring  $(A,\fm)$ (with the $\fm$-adic topology),
we have $\cX \cong \Spf A$.
Let
\[
 \cY = \lim _{\longrightarrow} \Spec A /I_n
\]
be a pseudo-subscheme where 
\[
A \supseteq I_1 \supseteq I_2 \supseteq \cdots
\]
is a descending chain of open ideals. 
Chevalley's theorem  \cite[Lem.\ 7]{Chevalley} (see also \cite[Ch.\ VIII, \S 5, Th.\ 13]{Zariski-Samuel-II})
says that either 
\begin{enumerate}
\item for every $n \in \NN$, there exists $i \in \NN$ with $I_i \subseteq \fm^n$, or
\item $\bigcap _i I_i \neq (0)$.
\end{enumerate}
In the former case, the $\{I_i\}$-topology coincides with the $\fm$-adic topology, and so
$\cY = \cX$.
In the latter case, replacing $A$ with $A/ \bigcap _i I_i$, we can reduce to the
former case. Consequently we see that $\cY = \Spf ( A/ \bigcap _i I_i ) $
and that $\cY$ is a subscheme of $\cX$.

We now prove the ``only if'' direction. 
Suppose that the underlying topological space of $\cX$ is not discrete.
Then there exists a closed but not open point $x$ of $\cX$.
Let $\Spf A \subseteq \cX$ be an affine neighborhood of $x$.
Then $\Spf A$ consists of at least two points. 
Let $A_\red$ be the reduced ring associated to $A$, 
that is, the ring $A$ modulo the ideal of nilpotent elements.
Then $\Spf A$ and $\Spf A_\red$ have the same underlying topologcial space. 
If $\hat A_\red $ is the $\fm $-adic completion of $A_\red$ with $\fm$ the maximal ideal
of $x$, then $\Spf \hat A_\red$ is a closed pseudo-subscheme of $\Spf A _\red$ consisting of a single point,
hence not isomorphic to $\Spf A_\red$.
Being injective, the natural map $A _\red \to \hat A_\red$ does not factors as $A_\red \to A_\red/ J \cong \hat A_\red $
 for any nonzero ideal $J$. Hence $\Spf \hat A_\red$ is not any closed subscheme of $\Spf A_\red$
 or of $\Spf A$.
\end{proof}

\subsection{The pseudo-closure of a pseudo-subscheme}

\begin{prop-defn}
Let $\cX$ be a  formal scheme and $\cY \subset \cX $ its pseudo-subscheme.
Suppose that the inclusion map of the underlying topological spaces
is quasi-compact.
Then there exists a smallest closed pseudo-subscheme $\bar \cY$ of $\cX$
that contains $\cY$ as an open subscheme. 
We call $\bar \cY$ the {\em pseudo-closure} of $\cY$ in $\cX$.
\end{prop-defn}

\begin{proof}
We first consider the case where $\cX$ is quasi-compact and $\cY$
 admits a basis of subschemes of definition,
say $\{Y_i\}_{i \in \NN}$. From
Proposition-Definition \ref{prop-closure-of-ordinary-subscheme}, for each $i$,
there exists the closure $\bar Y_i$ of $Y_i$ in $\cX$. 
Then we put 
\[
 \bar \cY := \varinjlim \bar Y_i .
\]
Let $\{Y'_i\}$ be another basis of subschemes of definition of $\cY$.
Then for every $i \in \NN$, there exists $j \in \NN$ such that
$\bar Y_i \subseteq \bar Y_j'$ and $\bar Y'_i \subseteq \bar Y_j$. 
Therefore $\bar Y_i'$ can be viewed as closed subschemes of $\bar \cY$
and also form a basis of subschemes of definition of $\bar \cY$.
It follows that $\bar \cY = \varinjlim \bar Y_i'$.
Thus $\bar \cY$ is independent of the choice of $\{Y_i\}$.
By construction, $\cY$ is an open subscheme of $\bar \cY$.
Moreover $\bar \cY$ is a smallest closed pseudo-subscheme of $\cX$
that contains $\cY$ as an open subscheme. 
Indeed if $\cZ$ is another closed pseudo-subscheme of $\cX$
containing $\cY$ as an open subscheme, then 
the $\bar Y_i$ are also closed ordinary subschemes of $\cZ$
and hence $\bar \cY$ is also a closed pseudo-subscheme of $\cZ$. 

We now consider the general case.
Then there exists an open covering $\cX = \bigcup_{\lambda \in \Lambda} \cX_\lambda$
such that for every $\lambda$, $\cX_\lambda$ is quasi-compact 
and  $\cY_\lambda := \cY \cap \cX_\lambda$ admits a basis of subschemes of definition.
For each $\lambda$, there exists the pseudo-closure $\bar \cY_\lambda$ of $\cY_\lambda$
in $\cX_\lambda$. Let $\cX _{\lambda \mu } := \cX_{\lambda} \cap \cX_{\mu}$,
 $\cY _{\lambda \mu} : = \cY \cap \cX _{\lambda \mu}  $ and $\bar \cY _{\lambda \mu}$
 the pseudo-closure of $\cY_{\lambda \mu}$ in $\cX_{\lambda \mu}$. 
Then by the construction above of $\bar \cY$,
$\bar \cY_{\lambda \mu} = \cY _\lambda \cap \cX_{\lambda \mu} $.
Therefore we can glue the $\bar \cY_\lambda$ and
obtain a closed pseudo subscheme $ \bar \cY $ of $\cX$
that contains $\cY$ as an open subscheme.
It is easy to see that $ \bar \cY$ is the smallest closed pseudo-subscheme with this property.
\end{proof}

\begin{expl}\label{expl-pseudo-closure}
Let $\cY$ be a  closed subscheme of $\Spf \CC [x^\pm,y][[t]]$
such that the only closed subscheme of $\Spf \CC[x,y][[t]]$ containing $\cY$
is $\Spf \CC[x,y][[t]]$. (Thanks to Theorem \ref{thm-example-no-closure},
such $\cY$ exists.) 
Then the pseudo-closure $\bar \cY$ of $\cY$ in $\Spf \CC[x,y][[t]]$
is not a closed subscheme of $\Spf \CC[x,y][[t]]$.
\end{expl}

\begin{rem}
Examples \ref{expl-infinite-embedded-points} and \ref{expl-pseudo-closure}
are both pseudo-subschemes that are not subschemes, but have different flavors.
It might be good to distinguish them, for example, by the following condition on a
 pseudo-subscheme $\cY \hookrightarrow \cX$: For any ordinary subscheme 
$Z \hookrightarrow \cX$, the fiber product $\cY \times _\cZ \cX$
is an ordinary scheme. While Example \ref{expl-infinite-embedded-points}
does not satisfy this, Example \ref{expl-pseudo-closure} does.
\end{rem}

\section{Formal separatrices of singular foliations}

In this section, we see that a pathological phenomenon of formal schemes
also comes from singularities of foliations.

\subsection{Formal separatrices}

Let $X$ be a smooth algebraic variety over $\CC$, and $ \Omega _{X}=\Omega_{X/\CC}$ 
the sheaf of (algebraic) K\"ahler differential forms.
A {\em (one-codimensional) foliation} on $X$ is an invertible saturated subsheaf
$\cF$ of $\Omega _{X}$ satisfying the integrability condition; $ \cF \wedge d\cF =0 $.
We say that a foliation $\cF$ is {\em smooth} at $x \in X$ if
the quotient sheaf $\Omega _X /\cF$ is locally free around $x$, and
that $\cF$ is {\em singular} at $x$ otherwise.
We say that $\cF$ is {\em smooth} if $\cF$ is smooth at every point.
The pair $(X,\cF)$ of a smooth variety  $X$ and a foliation on $X$
is called a {\em foliated variety}.
 
\begin{defn}
Let $(X,\cF)$ be a foliated variety, $x \in X (\CC)$, $ X_{/x} := \Spf \hat \cO_{X,x}$,
 $Y \subseteq  X_{/x}$  a closed subscheme of codimension one defined by $0 \neq f \in \hat \cO_{X,x}$,
and $ \omega \in \Omega_{X,x} $ a generator of $\cF_x$.
We say that $Y$ is  a {\em formal separatrix} (of $\cF$) at $x$ if
$f$ divides $ \omega \wedge d f $. 
\end{defn}

Because of Leibniz rule, 
 $Y$ is a formal separatrix if and only if its associated reduced formal scheme $Y_\red$ 
 is a formal separatrix. 
Frobenius theorem says that if $\cF$ is smooth at $x$,
there exists a unique smooth formal separatrix of $\cF$ at $x$.
Miyaoka \cite{Miyaoka} proved  that  the family of smooth formal separatrices
at smooth points of a foliation form a formal scheme:

\begin{thm}\label{thm-Miyaoka}\cite[Cor.\ 6.4]{Miyaoka}
Let $(X,\cF)$ be a foliated variety. Suppose that $\cF$ is smooth.
Then there exists a closed subscheme $\cL$ of  $(X \times _\CC X)_{/ \Delta_X}$ such that
for every point $x \in X$, $ p_2( p_1^{-1}(x)) $ is the smooth formal separatrix of $\cF$ at $x$.
Here $\Delta_X \subseteq X \times _\CC X$ is the diagonal and 
$(X \times _\CC X)_{/ \Delta_X}$ is the completion of $X \times _\CC X$ along $\Delta_X$.
\end{thm}

Let $(X,\cF)$ be a foliated variety and  $C \subseteq X$ a closed smooth subvariety of dimension 1.
Suppose that $C$ meets only at a single point $ o $ with the singular
locus of $\cF$. Let $U \subseteq X$ be the smooth locus of $\cF$
and $\cL \subseteq (U \times _\CC U)_{\Delta _U}$ the family of formal separatrices
as in the theorem.
Then $C \setminus \{o\}$ is a closed subvariety of $U$.
The fiber product 
\[
\cL_{C\setminus \{o\}} := ( C \setminus \{o\}) \times _{U, \pr_1} \cL
\]
 is the family of
the smooth formal separatrices over $C \setminus \{o\}$,
and a subscheme of $(C \times _\CC X)_{/\Delta_ C}$.
Let $\cL _C := \overline{ \cL_{ C\setminus \{o\} } }$ be the pseudo-closure of
$ \cL_{ C\setminus \{o\} }$ in $(C \times _\CC X)_{/\Delta_ C}$.

\begin{prop}\label{prop-equivalent-conditions-closed}
The following are equivalent:
\begin{enumerate}
\item ${\cL_{C}}$ is Noetherian.
\item ${\cL_{C}}$ is adic.
\item $ {\cL_{C}}$ is pre-Noetherian.
\item $ {\cL_{C}}$ is a closed subscheme of $( C \times _\CC  X ) _{/ \Delta_ C}   $.
\item ${\cL_{C}}$ is a closed pre-subscheme of $( C \times _\CC   X ) _{/\Delta_ C}   $.
\end{enumerate}
\end{prop}

\begin{proof}
$1 \Rightarrow  2 $  and $1 \Rightarrow 3$: Trivial.

$2 \Rightarrow 1$: It follows from Lemma \ref{lem-adic+pro-Noetherian-Noetherian}. 

$3 \Rightarrow 2$: 
The underlying topological space of $\cL_C$ is identified with that of $C$.
Shrinking $C$, we may suppose that $\cL_C$ is affine, say 
$\cL_C = \Spf A$ with $A$ a Noetherian admissible ring.
Let $ I \subseteq A $ be the largest ideal of definition.
This is a prime ideal and the symbolic powers $ I^{(n)} $
form a basis of ideals of definition in $A$. 

If $  f,g \in A $ are nonzero elements, then 
for $ n \gg 0$, their images $\bar f, \bar g $ in $A/I^{(n)}$
are nonzero. Their restrictions  $\bar f |_{C \setminus \{o\}}$ and $\bar g|_{C \setminus \{o\}}$ to $C \setminus \{o\}$
are also nonzero. Therefore  the restrictions $f |_{C \setminus \{o\}}$ and
 $ g|_{C \setminus \{o\}}$ of $f $ and $g$ are nonzero.
Consequently the restriction $ (fg) |_{C \setminus \{o\}}$
of the product $fg$ does not vanish, and the product $fg$ does not neither.
Thus $A$ is a domain. 

Let $\fm \subseteq A$ be the maximal ideal of $o$ and 
  $\hat A$ the $\fm$-adic completion of $A$.
We claim that $\hat A$ is also a domain.
To see this, we may suppose that $( C \times _\CC X )_{/\Delta _C}$ is affine,
say $\Spf B$. Let $\fn \subseteq B$ be the maximal ideal of $o$
and $\hat B$ the $\fn$-adic completion of $B$.
Put 
\[
 J_n := \Ker (\hat B \to \hat A / \hat I^{(n)}).
\]
Then  $\hat A \cong \varprojlim \hat B / J_n$.
Since $C$ is smooth, in particular, analytically irreducible, $\hat B /J_1$ is a domain.
Now we can prove the claim in the same way as above.

From \cite[page 33, Lem.\ 3]{Zariski} (see also \cite[Ch.\ VIII, \S 5, Cor.\ 5]{Zariski-Samuel-II}), 
the $I^{(n)}$-topology on $A$ coincides with the $I$-adic topology. 

$2 \Rightarrow 4 $: It is a direct consequence of Proposition \ref{prop-surjective-adic}.

$4 \Rightarrow 5$: Trivial.

$5 \Rightarrow 2$: Put $ \cX :=  ( C \times _\CC   X ) _{/\Delta_ C}   $ and $ \cY := \cL_C $.
Identifying the underlying topological spaces of $\cX$ and $\cY$ with $C$,
the point $o \in C$ can be viewed as a point of $\cX$ and $\cY$.
Since $\cX$ is Noetherian, from \cite[0, Cor.\ 7.6.18]{EGA}, the stalk $\cO_{\cX,o}$ is Noetherian.
Since tha natural map $\cO_{\cX ,o} \to \cO_{\cY,o}$ is surjective,
$\cO_{\cY,o}$ is also Noetherian.
Let $I \subseteq \cO_{\cY,o}$ be the ideal of the topologically nilpotent elements.
Then the symbolic powers $I^{(n)}$ of $I$ form a basis of open neighborhoods of
$\cO_{\cY,o}$. (Note that $\cO_{\cY,o}$ is not a priori complete.)
As in the proof of ``$3 \Rightarrow 2$'', we can show that the topology on $\cO_{\cY,o}$ is 
identical to the 
$I$-adic topology. 
Hence if $\cI \subseteq \cO_\cY$ is the ideal sheaf of topologically nilpotent sections,
then its powers $\cI^n$ form a basis of ideals of definition,
and $\cY$ is adic.
\end{proof}

When $\cL_C$ is a closed subscheme,
it allows us to take the limit of smooth formal separatrices along $C$:

\begin{thm}\label{thm-limit-formal-separatrix}
Suppose that one of the conditions in Proposition \ref{prop-equivalent-conditions-closed} holds.
Then the fiber of $  {\cL_C} \to C $ over $o$ is
a formal separatrix at $o$. 
\end{thm}

\begin{proof}
We need to use complete modules of differentials of
Noetherian formal schemes. 
For a morphism $f:\cY \to \cX$ of Noetherian formal schemes,
we have a complete module of differentials, $\hat \Omega_{\cY/\cX}$, 
which is a quasi-coherent $\cO_\cY$-module, and have a derivation 
$\hat d_{\cY/\cX}: \cO_\cY \to \hat \Omega_{\cY/\cX}$.
We refer to \cite{AJP} for details.

If necessary, shrinking $X$, we can take a nowhere vanishing 
$\omega \in \cF (X)$.
Let 
\[
 \psi:  \cX := (C \times _\CC X)_{/ \Delta_C} \to X.
\]
be the projection.
Pulling back $\omega$, we obtain a global section $\psi^* \omega$ of $ \hat \Omega_{\cX /C} $.
Since $\cL_C$ is a hypersurface in $\cX$, it is defined by a section $f$ of $\cO_\cX$. 
Since the restriction of $\cL_C$ to $C \setminus \{o\}$ is
 the family of formal separatrices along $C \setminus \{o\}$, $f$
 divides $\psi^* \omega \wedge \hat d_{\cX/C}f $. 

Let $ Y$ be the fiber of $\cL _C \to C$ over $o$, which is a hypersurface of $X_{/o}$
defined by the image $\bar f \in \hat \cO_{X,o}$ of $f$.
Then $\bar f$ divides $\omega \wedge \hat d_{X_{/o}/\CC} \bar f$. 
Hence $Y$ is a formal separatrix. 
\end{proof}

\subsection{Jouanolou's theorem and its application}

We recall Jouanolou's result on Pfaff forms.
We refer to  \cite{Jouanolou} for details.

An {\em  algebraic Pfaff form of degree $m$}  on $\PP_\CC ^2$ is a one-form
\[
 \omega = \omega _1 dx + \omega_3 dy + \omega_3 dz
\]
such that $\omega _i$ are homogeneous polynomials of degree $m$ and the equation
\[
  x \omega _1 + y \omega _2 + z \omega _3 =0
\] 
holds. 
A {\em Pfaff equation of degree $m$} on $\PP_\CC^2$ is a class of algebraic Pfaff forms 
modulo nonzero scalar multiplications. 

Let $\omega$ be an algebraic Pfaff form on $\PP_\CC^2$ and $ [\omega] $ its Pfaff equation class.
An {\em algebraic solution} of  $\omega$ or $[\omega]$ is a class of  homogeneous polynomials 
$f \in \CC[x,y,z]$
modulo nonzero scalar multiplications
such that $f$ divides $\omega \wedge df $.

Let $V_m$ be the vector space of the algebraic Pfaff forms of degree $m$ on $\PP_\CC^2$.
Then the set of the Pfaff equations of degree $m$ on $\PP_\CC^2$
is identified with the projective space $\PP (V_m)= (V_m \setminus \{0\} )/ \CC^*$.
Define 
\[
 Z_m \subseteq \PP (V_m)
\]
to be the set of the Pfaff equations that have no algebraic solution.

\begin{thm}\cite[\S 4]{Jouanolou}
Suppose $m \ge 3$. 
Then $Z_m$ is the intersection of countably many nonempty Zariski open subsets of $\PP(V_m)$
and contains the class of the algebraic Pfaff form
\[
 (x^{m-1}z -y^m ) dx + (y^{m-1}x -z^m) dy +(z^{m-1}y -x^m ) dz .
\]
\end{thm}

From \cite[page 4, Prop.\ 1.4]{Jouanolou}, every algebraic Pfaff form $\omega$ on $\PP_\CC^2$ is
integrable; $  d \omega \wedge \omega = 0$. So $\omega$ defines also a foliation $\cF_\omega$ on $\CC^3$.
From \cite[page 85, Prop.\ 2.1]{Jouanolou}, the only singular point of $\cF _\omega$
is the origin. Accordingly we can define 
the family $\cL_{\omega, C \setminus \{o\}}$ of formal separatrices along $C \setminus \{o\}$
and its pseudo-closure
$\cL _{\omega,C}$ for any line $C \subset \CC^3$ through the origin.

Let $ f = \sum_{i \ge n} f_i  \in \CC [[x,y,z]] $. Here $f_i$ is a homogeneous polynomial of degree $i$
and $f_n \ne 0$. Suppose that $f$ defines a formal separatrix at the origin, equivalently
that $f$ divides $\omega \wedge df$.
Then the class of $f_n$ is an algebraic solution of the Pfaff equation $[\omega]$.
Hence if $[\omega ] \in Z_m$, then $\cF_\omega$ has no formal separatrix at the origin.

\begin{cor}\label{cor-last-corollary}
For $[\omega] \in Z_m$ and a line $C \subseteq \CC^3$ through the origin, 
a pseudo-subscheme $\cL _{\omega , C}$ of $ (C \times _\CC \CC^3)_{\Delta _C} \cong \Spf \CC[w][[x,y,z]]$
is neither a closed pre-subscheme, pre-Noetherian nor adic.
\end{cor}

\begin{proof}
If $\cL_{\omega,C}$ is either a closed pre-subscheme, pre-Noetherian or adic,
then from Theorem \ref{thm-limit-formal-separatrix},
the foliation $\cF_\omega$ has a formal separatrix at the origin.
Hence $[\omega] \notin Z_m$.
\end{proof}



\end{document}